\documentclass [oneside]{amsart}
\usepackage {amscd}
\def\bc{\begin{center}}
\def\ec{\end{center}}
\def\px{\frac{\partial}{\partial x}}
\def\py{\frac{\partial}{\partial y}}
\def\pz{\frac{\partial}{\partial z}}

\def\EC{\mathcal{E}}
\def\EI{\mathcal{I}}

\newcommand{\C}{{\mathbb C}}
\newcommand{\CP}{{\mathbb {P}_{\mathbb{C}}^2}}
\newcommand{\CPP}{{\mathbb {P}_{\mathbb{C}}^{k-1}}}
\newcommand{\RP}{{\mathbb {R}P(2)}}
\newcommand{\CPu}{{\mathbb {C}P(1)}}
\newcommand{\PC}{{\mathbb {P}_{\mathbb{C}}^2}}
\newcommand{\CPPP}{{\mathbb {P}_{\mathbb{C}}^n}}
\newcommand{\cqd}{\ \hfill\rule[-1mm]{2mm}{3.2mm}}
\theoremstyle{change}
\newtheorem{theorem}{Theorem}
\newtheorem{proposition}{Proposition}

\newtheorem{definition}{Definition}
\newtheorem{remark}{Remark}
\newtheorem{corollary}{Corollary}
\newtheorem{example}{Example}
{ 
  \theoremstyle{margin}{
   }
}

\begin{document}

\title{Vector Fields, Invariant Varieties and Linear Systems}
\author{Jorge Vit\'orio Pereira}

\maketitle
\bc
Instituto de Matem\'{a}tica Pura e Aplicada, IMPA, Estrada Dona Castorina, 110 \\
Jardim Bot\^{a}nico, 22460-320 - Rio de Janeiro, RJ, Brasil. 
email : jvp@impa.br 
\ec

\begin{abstract}
We investigate the interplay between invariant varieties of vector fields and the inflection locus of linear systems with respect to the vector field. Among the consequences of such investigation we obtain a computational criteria for the existence of rational first integrals of a given degree, bounds for the number of first integrals on families of vector fields and a generalization of Darboux's criteria in the spirit of \cite{L}. We also provide a new proof of Gomez--Mont's result on foliations with all leaves algebraic, see \cite{GM}.   
\end{abstract}
%-------------------------------------------------------------------------------------------
% Introducao
%-------------------------------------------------------------------------------------------
\section{Introduction}
\rm
When studying algebraic curves one of the most fruitfull concept is the one of inflection and higher order inflection points. For a smooth plane curve, i.e. a smooth compact Riemann surface embedded in $\PC$, the inflection points are precisely the points were the tangent line has contact of order at least $3$ with the curve. If $C$ is a curve then a point $p \in C$ is  an  inflection point of order $d$  if there exist an curve of degree $d$ that has higher order contact with the curve $C$ at $p$. Here higher order contact means that the order of contact is at least the dimension of the vector space of polynomials of degree at most $d$. 

For plane curves, the inflections points are computed through the Hessian and the higher order inflection points are not so easy to obtain. As far as the author knows, the first mathematician to pursue the question of determining higher order inflection points for plane curves was Cayley, see \cite{Ca}. He succeed in giving a formula for the inflection points of order two. Although, the formula obtained by Cayley is not very simple.

According to Cukierman the problem of giving formulas for inflection points of order greater than three of plane curves does not seem to have been  solved in the classical literature. In \cite{Cu}, he gives an approach to obtain "almost explicit" formulas for the higher order inflection points of plane curves and complete intersection curves on some projective space.

The goal of the first part of this paper is to introduce and show how to compute inflection and higher order inflection points for holomorphic vector fields on the complex projective plane. In more concrete terms, given a vector field $X$ on $\PC$, we define effective divisors $\EC_d(X)$ on $\PC$, such that the restriction of the divisor to  any solution of the vector field $X$ coincides with the inflection points of order $d$ of the solution. 

In contrast with the case of plane curves, the formulas obtained for $\EC_d(X)$ are not very complicated. At first sight this seems to be paradoxical, but if $C$ is a smooth algebraic curve invariant by the algebraic vector field $X$, then the divisor $\EC_d(X)$ restricted to the curve gives something more than the inflection points of order $d$ of the curve $C$. The fact is that besides the inflection points of order $d$, $\EC_d(X)_{|C}$ also contains the singularities of the vector field along $C$.  
 
The initial motivation for introducing such concepts was to have a tool for detecting invariant algebraic curves of a given degree $d$ of a vector field $X$, and bound their number in function of $d$ and the degree of $X$. In fact using the divisors $\EC_d(X)$, which we call {\it extactic curves}, we obtain such bounds. These bounds turn out to be of different nature of Jouanolou's bound which are obtained through Darboux's approach, see \cite{J}.

Among the applications of the extactic curves, one can find a computational criteria for the existence of first integral of a given degree and some properties of families of foliations on $\PC$.  

Since we believe that the concepts here introduced may be useful for studying concrete examples of real and complex algebraic vector fields, we try to be as explicit as possible in the first part of paper.

In the second part we use a more intrinsic approach and generalize the concepts and some of the results for vector fields on arbitrary complex manifolds. 

Developing the concepts in such generality we show that we can use our methods to detect a class of solutions more general than the algebraic class. One illustrative result of the method is a generalization of Darboux's criteria for the existence of first integrals in the spirit of \cite{L}. Using the same sort of ideas we also obtain a new proof of a result proved in \cite{GM} by Gomez--Mont about foliations with all leaves algebraic.

I would like to thank P. Nogueira, for her help in working out the formal definition of extactic divisors, L.G. Mendes, for showing me  reference \cite{Fi}, and specially E. Esteves for many helpful discussions about linear systems, inflection points and algebraic vector fields.  

%-------------------------------------------------------------------------------------------
% Nocoes Basicas
%-------------------------------------------------------------------------------------------
%-------------------------------------------------------------------------------------------
% Curvas Invariantes
%-------------------------------------------------------------------------------------------
\part{Extactic Curves on the Projective Plane}\label{explicita}

\section{Affine and Projective polynomial Vector Fields}

If $X$ is a polynomial vector field on $\C^2$, then $X$ can be written, in a  unique way, as
\[
  X = a(x,y) \px + b(x,y) \py + g(x,y) \left( x \px + y \py \right) \, ,
\]
where $g$ is a homogeneous polynomial of degree $d$ and $a, b$ are polynomials of degree at most $d$. We define $d$ as the $degree$ of the vector field. When $g$ is identically zero we say that the line at infinity is invariant. 

If we consider the homogeneous polynomials 
\[
  \overline{a}(x,y,z) = z^d \cdot a \left( \frac{x}{z}, \frac{y}{z} \right) \, \, \, {\rm and } \, \, \, \,
\overline{b}(x,y,z) = z^d \cdot b \left( \frac{x}{z}, \frac{y}{z} \right) \, ,
\]
then the vector field 
\[
  \overline{X} = \overline{a}(x,y,z) \px + \overline{b}(x,y,z) \py + g(x,y,z) \pz
\]
is a projectivization of $X$. If $Y$ is any other projectivization of $X$ then $\overline{X} - Y$ is a multiple of the radial(or Euler) vector field $R = x \px + y \py + z\pz$. 

Reciprocally, if $Y$  is a polynomial vector field on $\C^3$ with homogeneous coefficients then $Y$ induces via the radial projection a field of directions on $\CP$. Observe that any homogeneous vector field on $\C^3$ whose difference with $Y$ is a  multiple of the radial vector field, induces the same field of directions on $\CP$. We abuse the language and say that such a field of directions is a projective vector field, or just a vector field on $\CP$.

\section{Higher order inflection Curves of Vector Fields on $\CP$}\label{ssimples}
\rm
In this section we define the  extactic curves, $\mathcal{E}_d$ for $d \in \mathbb{N}$, for vector fields on $\CP$. 
These curves describe the inflection and higher order inflection points for solutions of the vector field.

%-------------------------------------------------------------------------------------------
\subsection{Extactic  points  of plane curves}
\begin{definition}
\rm
A $n$-inflection point of a curve in $\mathbb{P}_{\mathbb{C}}^2$  is a point where 
the multiplicity of  intersection of the curve with some algebraic curve of 
degree $n$ is greater than
$$
  d(n) = \frac{n(n+3)}{2}.
$$
Note that $d(n)$ is the dimension of the space of plane curves of degree n.
\end{definition}

We use the term extactic point following  V.I. Arnold, see \cite{Ar}. 

\begin{remark} \label{comentario}
Observe that when $C$ is an algebraic curve in 
$\mathbb{P}_{\mathbb{C}}^2$ or $\mathbb{P}_{\mathbb{R}}^2$, if 
every point of $C$ is a $n$-inflection point then the degree of $C$ is at most $n$.
\end{remark}
%-------------------------------------------------------------------------------------------
\subsection{Extactic curves of vector fields}\label{extactic} \rm
Our purpose is to describe in an unified way the inflections and higher--order inflections points of the solutions of the vector field. If $X = a(x,y) \frac{\partial}{\partial x} + b(x,y)\frac{\partial}{\partial y}$ is a vector field on $\mathbb{C}^2$ and $(x,y(x))$ a parametrization of a solution we have that 
$ \det{
(\smallmatrix
  a & b  \\  1 & \frac{dy}{dx} 
\endsmallmatrix)
} = a \frac{dy}{dx}(x) - b = 0 $. Therefore, $ \frac{dy}{dx} = \frac{b}{a}$. 
To obtain the inflections of the curve $(x,y(x))$ we have to calculate the determinant: 
$$
 \det 
   \bmatrix
     1 & \frac{dy}{dx}  \\  0 & \frac{d^2y}{dx^2} 
   \endbmatrix
  = \frac{d}{dx} \frac{b(x,y(x))}{a(x,y(x))} = 
  \frac { (\frac{\partial a}{\partial x} + \frac{\partial a}{\partial y}\frac{a}{b} )b 
  - (\frac{\partial b}{\partial x} + \frac{\partial b}{\partial y}\frac{a}{b}  )a }{a^2}
$$

Hence, in the open set $\C^2 \setminus (a(x,y) = 0)$, the inflections points of the solutions 
describe a curve given by the expression above. In an analogous manner we can get equations 
for the inflection curve in the open set $\mathbb{C}^2 - (b(x,y) = 0)$.

To calculate the $2$-inflection points of the solutions we consider the image of the curve 
$(x,y(x))$ under the  $2$-Veronese map. Such image is the curve parametrized by 
$$
  (x,x^2,x \cdot y(x),y(x),y(x)^2)
$$ 
in $\C^5$. If we calculate its flattening points we obtain the $2$-inflection points of the solution. Such calculation can be done evaluating the determinant of the matrix,
$$
\bmatrix 
  1 & 2x & \frac{d xy(x)}{dx}    & \frac{d y(x)}{dx} &\frac{d y(x)^2}{dx} \\
  0 & 2  & \frac{d^2 xy(x)}{dx^2}& \frac{d^2 y(x)}{dx^2} &\frac{d^2 y(x)^2}{dx^2} \\
  0 & 0  & \frac{d^3 xy(x)}{dx^3}& \frac{d^3 y(x)}{dx^3} &\frac{d^3 y(x)^2}{dx^3} \\
  0 & 0  & \frac{d^4 xy(x)}{dx^4}& \frac{d^4 y(x)}{dx^4} &\frac{d^4 y(x)^2}{dx^4} \\
  0 & 0  & \frac{d^5 xy(x)}{dx^5}& \frac{d^5 y(x)}{dx^5} &\frac{d^5 y(x)^2}{dx^5} \\
\endbmatrix \, .
$$
The $n$-inflection points can be obtained through the use of the $n$-Veronese map in a completely similar way.

Until now we have been working locally, although it is possible to give global expressions, on $\CP$ for these n-inflection curves. 
If $X = A \frac {\partial}{\partial x} + B \frac {\partial}{\partial y} + C \frac {\partial}{\partial z}$ 
is a homogeneous vector field in $\mathbb{C}^3$ then the equation for the inflection curve, 
or the first extactic curve (which we will denote by $\mathcal{E}_1(X)$) 
of the induced foliation on $\mathbb{P}^2_{\mathbb{C}}$ is:
\begin{equation}\label{form1}
\mathcal{E}_1(X) = \det
\bmatrix 
  x & y & z \\
  X(x) & X(y)  & X(z) \\
  X^2(x) & X^2(y)  & X^2(z) 
\endbmatrix \, ,
\end{equation}
where $X^k(f) = X(X^{k-1}(f))$, for any polynomial $f$.
\begin{example}
\rm
Let  
\begin{eqnarray*}
   X = (x^3 - z^3)x \frac{\partial}{\partial x}  + (y^3 -z^3)y \frac{\partial}{\partial y} \ , \\ 
   Y = - y^2z^2 \frac{\partial}{\partial x} + -x^2z^2 \frac{\partial}{\partial y} 
       + x^2y^2 \frac{\partial}{\partial z}
\end{eqnarray*}
and $Z = tX + sY$, $(s,t) \in \C^2$ and $s\cdot t \neq 0$,  be a projectivization of Lins Neto's example \cite{LN}. 
Then  
$$
  \mathcal{E}_1(Z)=2 L_9 \cdot(t^2sy^3-xyzs^3+2xyzt^3+z^3t^2s+st^2x^3),
$$
where $L_9=(x^3 - y^3)(x^3 - z^3)(y^3 - z^3)$. Observe that the nine invariant lines 
for any $Z$  are contained in the first extactic curve of $Z$.
\end{example}

To understand why formula (\ref{form1}) works, suppose $p \in \C^3$ is a non--singular point of $X$. Here non-singular means that the vector field $X$ is not colinear with the radial vector at $p$, or in other terms $p$ is a non--singular point of the codimension one foliation of $\C^3$  induced by $X$ and $R$. By the existence of local solutions for ordinary differential equations there is a germ of curve $V$ around  $p \in \C^3$ which is a local orbit of $X$. As consequence the vector field $X$ restricts to $V$,which means that $X$ acts as a derivation on the local functions of $V$.  Since $V$,seen as a germ of projective curve, has dimension $1$ the restriction of $X$ to $V$ can be seen as the derivative of a local parameter $t$. Hence formula (\ref{form1}) can be written on $V$ as:
$$
\mathcal{E}_1(X)_{\mid V} (t) = \det
\bmatrix 
  x(t) & y(t) & z(t) \\
  \frac{\partial x(t)}{\partial t} &  \frac{\partial y(t)}{\partial t}  & \frac{\partial z(t)}{\partial t} \\
  \frac{\partial^2 x(t)}{\partial t^2} &  \frac{\partial^2 y(t)}{\partial t^2}  & \frac{\partial^2 z(t)}{\partial t^2} 
\endbmatrix \, ,
$$
and in fact represents the infection points of $V$ in a neighborhood of $p$.

Similarly, we have a global equation for the curve of 2-inflection points of the vector field, 
or the second extactic curve $\mathcal{E}_2(X)$, which is the determinant of :
$$
\bmatrix 
  x^2 & xy & xz & y^2 & yz & z^2 \\
  X(x^2) & X(xy) & X(xz) & X(y^2) & X(yz) & X(z^2) \\
  X^2(x^2) & X^2(xy) & X^2(xz) & X^2(y^2) & X^2(yz) & X^2(z^2) \\
  X^3(x^2) & X^3(xy) & X^3(xz) & X^3(y^2) & X^3(yz) & X^3(z^2) \\
  X^4(x^2) & X^4(xy) & X^4(xz) & X^4(y^2) & X^4(yz) & X^4(z^2) \\
  X^5(x^2) & X^5(xy) & X^5(xz) & X^5(y^2) & X^5(yz) & X^5(z^2) \\
\endbmatrix
$$ 

The $d$-th extactic curve $\mathcal{E}_d(X)$ can be described in a completely similar way. The equation of $\mathcal{E}_d(X)$ is given by the determinant of 
a matrix such that the first row is formed by a basis of the monomials in $x,y$ 
and $z$ of degree $d$, and the $i$-th row is the derivation $X$ applied in the $(i-1)$-th row.

\begin{proposition}\label{pini}
Every  algebraic curve of degree $n$  invariant by the vector field $X$ is a factor of $\EC_n(X)$.
\end{proposition}

\noindent{\it Proof}: Let $F$ be an invariant algebraic curve of degree $n$.
Since  the choice of the basis of the $\C$--vector space plays no role in the definition
of extactic curve,  we can choose a basis where $F$ appears. Since 
\[
\begin{array}{lclcl}
X(F)   &=& L_F F \ ,     & &                                     \\
X^2(F) &=& X(L_F F)      &=& \left( L_F^2 + X(L_F) \right) F \ , \\
X^k(F) &=& X(X^{k-1}(F)) &=& \left ( \, \mbox{polynomial} \, \right) F         \ ,
\end{array}
\]
where $L_F$ is a polynomial, one can see that $F$ is a factor of $\EC_n(X)$. 
\cqd

\begin{theorem}\label{tini}
$X$ admits a first integral of degree $d$, but do not admit a first integral of degree smaller than $d$ if, and only if, $\mathcal{E}_d(X) = 0$ and $\mathcal{E}_{d-1}(X) \neq 0$.
\end{theorem}

\noindent{\it Proof}:
Let $p \in \mathbb{P}^2_{\mathbb{C}}$ be a non-singular point of $X$. 
Suppose that the solution passing through $p$ is parametrized, locally, by $(x,y(x))$. 
Since $\mathcal{E}_d(X)$ vanishes identically, the composition of our local solution with 
the $d$-Veronese map is contained in a hyperplane, so  $(x,y(x))$ must be contained in an 
algebraic curve of degree at most $d$. Since every leaf is algebraic it follows from Darboux's Criteria, see \cite{L}, that $X$ admits a first integral of degree at most $d$.
The fact that $\mathcal{E}_{d-1}(X) \neq 0$ implies  that the generical solution is of degree at least $d$. 

If $X$ admits a first integral of degree $d$ then every invariant curve is of degree at 
most $d$ and hence every point is a $d$-inflection point, i.e., $\mathcal{E}_d(X) = 0$. 
Since not every invariant curve has degree $d-1$,  $\mathcal{E}_d(X) \neq  0$.
\cqd

From the Theorem \ref{tini} and the Proposition \ref{pini} we derive:

\begin{theorem}\label{curvefinder}
Let $X$ a vector field in $\mathbb{P}_{\mathbb{C}}^{2}$. 
For every $d \in \mathbb{N}$ the equations of possibles invariant curves of degree less 
than or equal to $d$ appear as factors of $\mathcal{E}_d(X)$, and 
if $\mathcal{E}_d(X) = 0 $ then $X$ has a meromorphic first integral of degree at most $d$.
\end{theorem}

%-------------------------------------------------------------------------------------------
% Algoritmo
%-------------------------------------------------------------------------------------------
%-------------------------------------------------------------------------------------------
%Contagem
%-------------------------------------------------------------------------------------------

\section{Counting algebraic solutions and field of definition of the invariant curves}

Jouanolou in \cite{J} shows that a vector field of degree $d$ on $\CP$ without rational first integral has at most $d(d+2)/2$ irreducible algebraic solutions. Observe that he does not make any assumption on the degree of the algebraic leaves. 

Using the extactic curves we are able to obtain different bounds for the number of irreducible algebraic solutions. If $X$ is a vector field on $\CP$ and we denote by $n_i(X)$ the number of irreducible algebraic solutions of degree $i$ of $X$ the following proposition holds.

\begin{proposition}\label{cota}
Let $X$ be a homogeneous vector field, $\mathbb{P}_{\mathbb{C}}^{2}$, of degree $d$. 
If it does not have a first integral of degree $\leq n$ then 
$$
\sum_{i=1}^{n} i\cdot n_i(X) \leq \frac{d(n^4+6n^3+11n^2+6n)-n^4-2n^3+n^2+2n}{8}.
$$
\end{proposition}

\noindent{\it Proof}:
One has just to observe that the left hand side of the inequality is bounded by the degree of
$\mathcal{E}_n(X)$. 
Since the space of monomials in $x,y,z$ of degree $n$ is $\pmatrix n+2 \\ 2 \endpmatrix$, we have
that the degree of $\mathcal{E}_n(X)$ is  
$$
\sum_{i=0}^{ (\smallmatrix n+2 \\ 2 \endsmallmatrix) -1 }{i(d-1) + n} 
= ((\smallmatrix n+2 \\ 2 \endsmallmatrix))n + 
(d-1) \frac{(\smallmatrix n+2 \\ 2 \endsmallmatrix)( (\smallmatrix n+2 \\ 2 \endsmallmatrix) -1 )}{2} = 
$$
$$
= \frac{d(n^4+6n^3+11n^2+6n)-n^4-2n^3+n^2+2n}{8}
.
$$
\cqd

\begin{corollary}
Let $X$ be a homogeneous vector field, $\mathbb{P}_{\mathbb{C}}^{2}$, of degree $d$. 
If it does not have a rational first integral of degree $\leq n$ then it has at most 
$$
\frac{d(n^3+6n^2+11n+6 )-n^3-2n^2+n+2}{8}
$$
invariant curves of degree $n$.
\end{corollary}
\noindent{\it Proof}:
If $C$ is an invariant curve of degree $n$ then $C$ is contained in 
$\mathcal{E}_n(X)$,  so we have at most 
$$\frac{deg(\mathcal{E}_n(X))}{n} 
=\frac{d(n^3+6n^2+11n+6 )-n^3-2n^2+n+2}{8}$$
invariant curves of degree $n$.\cqd

\begin{example}\rm
Applying Proposition \ref{cota} to bound the number of invariant lines of a vector field $X$, one can see that
\[
  n_1(X) \le 3\cdot {\rm deg}(X) \, .
\]
This bound turns out to be sharp.  For example, for each $d \in \mathbb{N}$, consider the vector fields $X_d$ on $\CP$, given in homogeneous coordinates by 
\[
  X_d = \left( x^{d-1} - z^{d-1} \right) x \px + \left( y^{d-1} - z^{d-1} \right) y \py \, .
\]
Then $X_d$ leaves invariant the algebraic curve $C_d$, cutted out by the polynomial 
\[ 
   F_d = xyz(x^{d-1} - z^{d-1})(y^{d-1} - z^{d-1})(x^{d-1} - y^{d-1}) \, .
\] 
Since $F_d$ can be expressed as a product of $3d$ distinct lines and $F_d = \EC_1(X_d)$, we conclude that the vector field $X_d$ admits exactly $3d$ invariant lines. 

Observe that the same bound holds for real vector fields on $\RP$. Although, the sharpness fails. To understand better why this happens see \cite{LL}.\qed
\end{example}

Using the same methods we can obtain smaller bounds if we try to count, for example, straight lines passing through a given point $p$. Instead of looking for the solutions in the $3$--dimensional vector space of lines, we have just to look for in the codimension one subspace of lines passing through $p$. For example, if $p = [ 0: 0: 1] \in \PC$ then equation of any invariant line passing through $p$  will appear as a factor of the determinant of the following matrix:
\[
\bmatrix 
  x & y  \\
  X(x) & X(y) 
\endbmatrix \, .
\]
Hence, if the $X$ has degree $d$, then the number of invariant lines passing through a given point $p$ is at most $d+1$.

Suppose now that we have an algebraic vector field on $\C^2$, and that its coefficients are in a normal subfield $\mathbb{K}$ of $\C$. The following question naturally arises: what can be said about the field of definition of the invariant algebraic curves? In other terms: what is the smallest extension $\mathbb{L}$ of $\mathbb{K}$ such that any invariant curve can be defined by a polynomial with coefficients in $\mathbb{L}$?

We can use the bounds for the number of invariant curves  to obtain bounds on the algebraic degree of the field of definition of algebraic invariant curves. 

\begin{proposition}
Let $\mathbb{K} \subset \C$ be a normal extension of $\mathbb{Q}$ and $X = a \px + b \py$ be a vector field on $\C^2$ of degree $d$, where $a, b \in \mathbb{K}[x,y]$. Suppose that exists an invariant algebraic curve of degree $l$, cutted out by a polynomial $f \in \mathbb{L}[x,y]$, where $\mathbb{L}$ is a normal extension of $\mathbb{K}$ contained in $\mathbb{C}$, that cannot be defined in any normal subfield of $\mathbb{L}$. If 
\[
  [\mathbb{L}:\mathbb{K}] > {\rm min}\left(\frac{d(n^4+6n^3+11n^2+6n)-n^4-2n^3+n^2+2n}{8},                                       \frac{d(d+2)}{2}\right) \, ,
\] 
then $X$ admits a first integral.
\end{proposition}

\noindent{\it Proof}:
Since $f$ is invariant by $X$, we have that $X(f)=L_f \cdot f $, for some polynomial $L_f \in \mathbb{L}[x,y]$. Applying the Galois automorphisms of the extension $[\mathbb{L}:\mathbb{K}]$, we obtain $[\mathbb{L}:\mathbb{K}]$ distinct invariant algebraic curves of degree $l$. Hence the Theorem follows by Jouanolou's bound and Proposition \ref{cota}.
\cqd

%-------------------------------------------------------------------------------------------
% Familias de Folheacoes
%-------------------------------------------------------------------------------------------

\section{Families of holomorphic foliations on $ \PC $ }

In \cite{LN}, Lins--Neto shows the existence of some very special families of foliations on $\PC$ parametrized by  the projective line. One interesting properties of such families is that not all foliations admits a rational first integral, but for a dense set $E$ in the parameter space the corresponding foliation has a rational first integral. Since the set in the parameter space admitting a rational first integral of degree at most $d$, for a fixed positive integer $d$, is algebraic, we have that the rational first integrals in the family have unbounded degree and $E$ admits the filtration 
\[
  E = \bigcup_{d \in \mathbb{N}}{E_d}\, ,
\]
where  $ p \in E_d$ if, and only if  the foliation corresponding to  $p$   admits a rational first integral of degree at most $d$.

In remark $4$ of \cite{LN}, Lins Neto says that "would be interesting to know what kinds of properties this set has". Here, we use the extactic curves to bound the growth of the cardinality of $E_d$ for any family of foliations parametrized by a projective line.

\begin{definition}\rm
Let $C \subset Fol(k)$ be a algebraic curve included in the space  of foliations of degree $k$.
We define the counting function of $C$, 
\[
  \pi_C : \mathbb{N} \rightarrow \mathbb{N} \cup +\infty,
\]
by the following rule : $\pi_C(d) = n$ if the number of points  in $C$ representing a 
foliation with rational first integral of degree at most $d$ is exactly $n$.
\end{definition}

\begin{example} \rm
Let $X = x \px$ and $Y= y \py$ be vector fields on $\C^2$. The family of vector fields $tX + sY$ can be seen
as $\CPu$ linearly embedded in $Fol(1)$, which we will denote by $C$. Whenever the ratio of $t$ and $s$ is a
rational number then $tX + sY$ admits a rational first integral. Suppose  that $t/s = p/q$ and $p$ and $q$ do
not have common factors. If $t/s$ is positive then the degree of the first integral is the maximum between $p$
and $q$, otherwise it is $|p|+|q|$. 

From the considerations of the previous paragraph we can show that 
\[
  \pi_C(n) \le K\cdot n^2.
\]
\end{example}

Considering a family of foliations parametrized by a projective line linearly embedded in the space of foliations, we obtain:

\begin{proposition}\label{simples}
Suppose $C$ is a $\CPu$ linearly embedded in $Fol(k)$. If $\pi_C(d) < \infty$ for every $d \in \mathbb{N}$, then  there exist a constant $K$ such that
\[
  \pi_C(d) \leq Kd ^4. 
\]
\end{proposition}

\noindent{\it Proof}: If we take two distinct projective vector fields in $C$, say $X$ and $Y$, we can recover $C$ by considering the linear combinations $sX + tY$. From the definition of $\EC_d$, one can see that:
\[
  \EC_d(sX + tY) = \sum_{   \alpha+\beta+\gamma={\rm deg} (\EC_d(X)}           
                       P_{\alpha,\beta,\gamma}(s,t)x^{\alpha}y^{\beta}z^{\gamma} \, .
\]
Since $\pi_C(d) < \infty$ then there exist a triple $(\alpha_0,\beta_0,\gamma_0)$ such that $P_{\alpha_0,\beta_0,\gamma_0}$ does not vanish identically. Hence the number of projective parameters  $(s:t)$ that has first integral of degree at most $d$ is bounded by the degree of $P_{\alpha_0,\beta_0,\gamma_0}$.

\cqd

%-------------------------------------------------------------------------------------------
% Secção 5
%-------------------------------------------------------------------------------------------

\part{Extactic Divisors on Complex Manifolds}

\section{Extactics divisors for holomorphic foliations }

Here we generalize some of the results of Part \ref{explicita} to foliations by curves on 
arbitrary non-singular complex manifolds. In order to do this, is imperative the reformulation of the concepts in a more intrinsic way, and to accomplish that we use freely the language of algebraic geometry(for example line bundles, tensor products and so on).

\subsection{Holomorphic foliations as morphisms}

\rm
Let $M$ be a complex manifold. An $1$-dimensional holomorphic foliation is 
given by the following data 
\begin{itemize}
  \item an open covering $\mathcal{U}={U_i}$ of $M$;
  \item  for each $U_i$ an holomorphic vector field $X_i$ ;
  \item for every non-empty intersection, $U_i\cap U_j \neq \emptyset $, a 
        holomorphic function $g_{ij} \in \mathcal{O}_M^*(U_i\cap U_j)$;
\end{itemize}
subject to the conditions :
\begin{itemize}
\item $X_i = g_{ij}X_j$ in $U_i\cap U_j$
\item $g_{ij}g_{jk} = g_{ik}$ in $U_i\cap U_j\cap U_k$.
\end{itemize}

If we denote by $\mathcal{L}$ the line bundle defined by the cocycle ${g_{ij}}$ we 
can understand the collection $X_i$ as a holomorphic section $\sigma$ of the bundle 
$TM \otimes \mathcal{L}$. Such section induces a morphism that goes from the cotangent 
bundle, denoted by $\Omega_M^1$, to the line bundle $\mathcal{L}$. This morphism, in the open set $U_i$, is  given by the interior product with the vector field $X_i$. 

Hence given an holomorphic foliation $\mathcal{F}$ we have an morphism:
$$
  \Phi_\mathcal{F} : \Omega_M^1 \rightarrow \mathcal{L}
$$

Reciprocally, given such morphism $\Phi$ we can canonically associate a holomorphic foliation to it. We leave 
the details for the reader. 

\subsection{Jet Bundles and Extactic Divisors}\label{jatos}
\rm
Suppose now that we have a foliation $\mathcal{F}$ on $M$, i.e.,  a morphism  $\Phi_\mathcal{F}$ of $\Omega^1_M$ to an invertible sheave $\mathcal{E}$. If we have a linear system $V \subseteq H^0(M,\mathcal{L})$ we are going to define {\it the extactic divisor, $\mathcal{E}(\mathcal{F},V)$, with respect to $V$}.  The extactic divisor $\mathcal{E}(\mathcal{F},V)$ can be understood geometrically as the inflection locus of the linear system $V$ with respect to the morphism $\Phi$.
 
First of all, consider the local Taylor expansion of a section $s \in V$ with respect 
to the vector field defining $\mathcal{F}$. Formally if you have 
a morphism 
$$
\Phi_\mathcal{F} : \Omega^1_M \rightarrow \mathcal{E}
$$
and a linear system $V \subseteq H^0(M,\mathcal{L})$, we choose a covering 
$\mathcal{U}$ of $M$ which trivializes both $\mathcal{L}$ and  $\mathcal{E}$. 
In a open set $U \in \mathcal{U}$ we can consider the morphism 
\begin{equation*}
  T^{(k)}_{|U} : H^0(M,\mathcal{L}) \otimes \mathcal{O}_{U} \rightarrow \mathcal{O}_{U}^k
\end{equation*}
defined by 
\begin{equation}\label{tlocal}
 T^{(k)}(s) = s + X_\mathcal{F}(s) \cdot t + X_\mathcal{F}^2(s) \cdot \frac{t^2}{2!} 
 + \cdots + X_\mathcal{F}^k(s)\cdot \frac{t^k}{k!}  
\end{equation}
where $X_\mathcal{F}(\cdot) = \Phi_\mathcal{F}(d(\cdot))$ and $s \in \mathcal{O}_{U}$ 
is an element of $H^0(M,\mathcal{L}) \otimes \mathcal{O}_{U}$ expressed in the chosen trivialization.

If we take open sets $U_\lambda \in \mathcal{U}$ we have that 
\begin{eqnarray*}
  \mathcal{L}_{|U_\lambda} &=& \mathcal{O}_{U_\lambda} \cdot \alpha_\lambda   \, , \\
  \mathcal{E}_{|U_\lambda} &=& \mathcal{O}_{U_\lambda} \cdot \beta_\lambda \, .
\end{eqnarray*}

Hence, for any $s_\lambda \in H^0(M,\mathcal{L}) \otimes \mathcal{O}_{U_\lambda}$, we obtain :
\begin{eqnarray*}
  s_\lambda &=& s_\lambda^{(0)} \cdot \alpha_\lambda , \\
  X_\mathcal{F}(s_\lambda) &=& X_\mathcal{F}(s_\lambda^{(0)}) \cdot \beta_\lambda 
     = s_\lambda^{(1)} \cdot \beta_\lambda \, .
\end{eqnarray*}

And generally,
$$
  X_\mathcal{F}^k(s_\lambda) = X_\mathcal{F}(s_\lambda^{(k-1)}) \cdot \beta_\lambda 
  = s_\lambda^{(k)} \cdot \beta_\lambda \, .
$$

When $U_\lambda \cap U_\mu \neq \emptyset$ :
\begin{eqnarray*}
  s_\lambda &=& s_\lambda^{(0)} \cdot \alpha_\lambda = l_{\lambda \mu}s_\mu^{(0)}
              \cdot \alpha_\mu \, , \\
  X_\mathcal{F}(s_\lambda) &=& s_\lambda^{(1)} \cdot \beta_\lambda 
                     =  X_\mathcal{F}(s_\lambda^{(0)} )\cdot \beta_\lambda =
         (X_\mathcal{F}(l_{\lambda \mu})\cdot s_\mu^{(0)} + l_{\lambda \mu}\cdot
            s_\mu^{(1)}
          ))\cdot e_{\lambda \mu} \cdot \beta_\mu \, .
\end{eqnarray*}

Or in matrix notation
$$
  \bmatrix s_\lambda^{(0)} \\ s_\lambda^{(1)} \endbmatrix = 
  \bmatrix 
    l_{\lambda \mu} & 0  \\
     X_\mathcal{F}(l_{\lambda \mu}) \cdot e_{\lambda \mu} & l_{\lambda \mu}   \cdot e_{\lambda \mu}  \\
  \endbmatrix
\cdot \bmatrix s_\mu^{(0)} \\ s_\mu^{(1)} \endbmatrix 
$$

By analogous computations one can show that :
\begin{equation}\label{transicao}
\bmatrix s_\lambda^{(0)} \\ s_\lambda^{(1)} \\ s_\lambda^{(2)} \\ \vdots \\ s_\lambda^{(k)}\endbmatrix = 
\bmatrix 
  l_{\lambda \mu} & 0 & 0 & 0 & 0  \\
   X_\mathcal{F}(l_{\lambda \mu}) \cdot e_{\lambda \mu} & l_{\lambda \mu} \cdot e_{\lambda \mu} & 0 & 0 & 0 \\
  \ddots & \ddots &  l_{\lambda \mu} \cdot e_{\lambda \mu}^2 & 0 & 0 \\
  \ddots & \ddots & \ddots & \ddots & 0  \\
  \ddots &  \ddots & \ddots & \ddots &  l_{\lambda \mu} \cdot e_{\lambda \mu}^k \\ 
\endbmatrix
\cdot \bmatrix s_\mu^{(0)} \\ s_\mu^{(1)} \\ s_\mu^{(2)} \\ \vdots \\ s_\mu^{(k)} \endbmatrix 
\end{equation}

Therefore, we define the vector bundle $J^k_{X_\mathcal{F}}\mathcal{L}$ as the vector bundle with rank $k+1$ and 
transition functions given by the matrix above. Now we are ready to consider the global {\it Taylor expansion} of the sections in $V$ with respect to $\mathcal{F}$. More precisely, we are able to patch together the morphisms of the form (\ref{tlocal}) to obtain just one morphism :
\begin{equation}\label{tglobal}
T^{(k)} : H^0(M,\mathcal{L}) \otimes \mathcal{O}_M \to J^k_{X_\mathcal{F}}\mathcal{L}
\end{equation}

We are almost ready to describe the inflection points of the linear system $V \subset H^0(M,\mathcal{L})$ 
with respect to $X_\mathcal{F}$. If we set $k$  as the dimension over $\mathbb{C}$ of $V$, and 
take the determinant of (\ref{tglobal}),
\begin{equation}\label{quase}
{\rm det } \,  T^{(k)} : \Lambda^{k} V \otimes \mathcal{O}_M \rightarrow  \Lambda^{k} J^{k-1}_{X_\mathcal{F}}\mathcal{L} \, .
\end{equation}
then after tensorizing we obtain a section of $\Lambda^{k} J^{k-1}_{X_\mathcal{F}}\mathcal{L} \otimes (\Lambda^{k} V)^{*}$
\begin{equation}\label{secao}
   \mathcal{O}_M \rightarrow \Lambda^{k} J^{k-1}_{X_\mathcal{F}}\mathcal{L} \otimes (\Lambda^{k} V)^{*} \ .
\end{equation}

Finally :

\begin{definition} \rm
The extactic curve of $\mathcal{F}$ with respect to the linear system $V$, 
$\mathcal{E}(\mathcal{F},V)$ is given by the zero locus of the section (\ref{secao}).
\end{definition}

%\begin{proposition}
%The following sequence is exact : 
%$$
%0 \rightarrow \mathcal{L} \otimes \mathcal{E}^{\otimes k} \rightarrow %J^k_{X_\mathcal{F}}\mathcal{L} \rightarrow J^{k-1}_{X_\mathcal{F}}\mathcal{L} \rightarrow 0 
%$$
%\end{proposition}
%\noindent{\it Proof}:
%Follows easily from the form of the transition matrix. \cqd

\begin{remark} \rm
The extactic curves, $\mathcal{E}_d(X)$, defined in subsection \ref{extactic} coincide with the extactics curves of the foliation induced by $X$ with respect to the linear system $H^0(\mathbb{P}^2_{\mathbb{C}}, \mathcal{O}_{\mathbb{P}^2_{\mathbb{C}}}(d))$. To verify this fact one has just to check locally, and is not very hard to show that locally the constructions are the same.
\end{remark}

\subsection{Order of Contact and Extactic Ideals}

So far we have defined the extactic divisors for $1$--dimensional holomorphic foliations on complex manifolds. One of the motivations to do that is to obtain an analogous of Theorem \ref{curvefinder}. In fact by a very similar argument we obtain the following:

\begin{proposition}\label{geral}
Let $\mathcal{F}$ be a $1$--dimensional foliation on the complex manifold $M$. If $V$ is a finite dimensional linear system, then every invariant hypersurface contained in the zero locus of some element of $V$, must be contained in the zero locus of $\EC(\mathcal{F},V)$. \end{proposition}

Since we are working on a manifold of arbitrary dimension, the extactic divisors are far from detecting precisely the invariant curves. To overcome this difficult we shall now introduce the extactic ideals. If we fix a finite dimensional linear system $V$, the main idea is to consider the points such that the contact with $V$ is of infinite order. 

Before defining the extactic ideals, let's make the notion of contact more precise. 

\begin{definition}\label{contato} \rm
Let $X_\mathcal{F}$ be an holomorphic vector field on the complex manifold $M$ and $s$ a holomorphic section of some line bundle. We say that the solution through $p$ has {\it contact of order $k$ with} $s$ when $k$ is the least non negative integer such that the radical of the ideal generated by $\overline{s}, \overline{X(s)}, \overline{X^2(s)}, \ldots, \overline{X^k(s)}$ is the local ring $\mathcal{O}_{M,p}$. If such a $k$ does not exist we say that $s$ has {\it flat contact} with $X$. Here $X$ means a local representant of $X_\mathcal{F}$ on a suitable open set and $\overline{s}$ means the image of any local representant of $s$ under the canonical morphism $\mathcal{O}_M \to \mathcal{O}_{M,p}$. We shall denote by $\nu(s,X_\mathcal{F},p)$ the contact of $s$ with $X_\mathcal{F}$ at the point $p$.
\end{definition}

Observe that any $s$ vanishing at a singular point $p$ of $X$, has flat contact with $X$ at $p$.

\begin{example} \rm
Let $X = \px$ be a polynomial vector field on $\C^3$ and $f$ a polynomial. If $f$ does not vanish at the origin then $\nu(f,X,0) = 0$, otherwise $\nu(f,X,0) \ge 1$. For example if 
\[ 
   f(x,y,z)= \sum_{i=k}^{n} f_i(y,z)x^i  \, ,
\]
with $f_k(0,0) \neq 0$, then $\nu(f,X,0) = k$.
\end{example}

\begin{proposition}
Let $\mathcal{F}$ be a $1$--dimensional foliation on the complex manifold $M$. Suppose $V$ is a finite dimensional linear system and $s \in V$. If $\nu(s,\mathcal{F},p) \ge {\rm dim}_\C V$ then $p \in \EC(\mathcal{F},V)$. Reciprocally if $p \in \EC(\mathcal{F},V)$ then exists an element $s \in V$ such that $\nu(s,\mathcal{F},p) \ge {\rm dim}_\C V$.
\end{proposition}
\noindent{\it Proof}: Choose a basis of $V$ starting with $s$. Using Lagrange's rule to expand $\EC(\mathcal{F},V)$ in terms of the first column of the matrix one can see that  $\EC(\mathcal{F},V)$ belongs to the maximal ideal correspondent to the point $p$.

When $\EC(\mathcal{F},V)(p)=0$, we have that the columns of the matrix used to compute $\EC(\mathcal{F},V)$ are linearly dependent in the point $P$. But that means that exists an element $s \in V$, such that $s(p) = X(s)(p) = \ldots = X^l(s)(p) = 0$, where $l=\dim_\C V$. Hence the result follows.
\cqd

Let $X_{\mathcal{F}}$ be an $1$--dimensional foliation on the manifold $M$. If $\mathcal{L}$ is an invertible sheaf then we can consider, as in \ref{jatos}, the morphisms 
\[
  T^{(l)} : H^0(M,\mathcal{L}) \otimes \mathcal{O}_M \to J^l_{X_{\mathcal{F}}}\mathcal{L} \, ,
\]
for any positive integer $l$. Suppose $V$ is a finite dimensional vector space contained in $H^0(M,\mathcal{L})$ and $k$ is its dimension.  Now, after taking the determinant and tensorizing by $(\Lambda^k V)^*$, we obtain sections
\[
  \sigma_l : \mathcal{O}_M \to \Lambda^k J^l_{X_\mathcal{F}} \mathcal{L} \otimes (\Lambda^k H^0(M,V)^* \, .
\]

\begin{definition}\label{ideal} \rm
The sheaf of ideals generated by $\ker \sigma_l$, where $l$ is any positive integer, is the 
extactic ideal of $X_{\mathcal{F}}$ with respect to the linear system $V$. We shall denote it by $\EI(X_{\mathcal{F}},V)$. 
\end{definition}

\begin{example}
If $X$ is a vector field on $\PC$ the  extactic ideal of $X$ with respect to the linear system $V=H^0(\PC,\mathcal{O}_{\PC}(k))$, $\EI(X,V)$, 
is generated, in homogeneous coordinates, by  $\sigma_{(k_1,\ldots,k_l)}$, where   
\begin{equation*}
\sigma_{(k_1,\ldots,k_l)} =
\mbox{det} \begin{pmatrix}  X^{k_1}(v_1)    &    X^{k_1}(v_2)       & \cdots & X^{k_1}(v_l)          \cr
                           X^{k_2}(v_1)       & X^{k_2}(v_2)       & \cdots & X^{k_2}(v_l)       \cr
                           \vdots       & \vdots       & \cdots & \vdots       \cr
                           X^{k_l}(v_1) & X^{k_l}(v_2) & \cdots & X^{k_l}(v_l) \end{pmatrix},
\end{equation*}
where  $0 \le k_1 < k_2 < \cdots < k_l$ and the $k_i$'s are integers,  $v_1, v_2, \ldots, v_l $ 
is a basis of  the space of homogeneous polynomial of degree $k$ in three variables.
\end{example}

\begin{proposition}\label{tideal}
The subvariety $N$ of $M$ associated to the radical of any extactic ideal of $X_{\mathcal{F}}$ is invariant.
\end{proposition}
\noindent{\it Proof}: Fix a linear system $V$ on $M$. If every element of $\EI(X,V)$ vanishes at $p$, then $p$ belongs to the radical ideal of $X_{\mathcal{F}}$ associated to $V$. Hence for every natural number $l \ge \dim_\C V$ one has a linear subspace $L_l \subset V$ such that any  nonzero element $f \in L_l$ satisfies $f(p) = X(f)(p) = \ldots = X^l(f)(p) = 0$. Being  the projectivization of $V$ compact we can find a element $f_{\infty}$ which vanishes at $p$ together with all its derivatives with respect to the local vector field $X$. In other words $f_{\infty}$ has flat contact with $X_\mathcal{F}$ at $p$.

Hence $N$ can be identified with the set of points $p \in M$ such that there exists a element of $V$ with flat contact  with $X_\mathcal{F}$ at $p$, and the invariance follows.

\cqd

\begin{example} \rm
In \cite{J}, Jouanolou proved that the vector fields on $\PC$ given in homogeneous coordinates by 
\[
  X_d = y^d \px + z^d \py + x^d \pz \, ,
\]
do not admit any invariant algebraic curve for any integer $d$ greater than $1$. 

To consider any linear system $V$ on $\PC$, is the same as to consider a finite dimensional space of homogeneous polynomials $V$ in three variables. Since we do not have any invariant algebraic leaves  for $X_d$ and the singular points of $X_d$ are clearly invariant, we have that the radical of the  extactic ideal $\EI(X_d,V)$ is exactly the ideal defining the singular set of $X_d$, for any $d \ge 2$. 
\end{example}

\begin{remark}\label{terminologia} \rm
If $\mathcal{F}$ is a holomorphic foliation on a complex manifold $M$, then, for us, a first integral for $\mathcal{F}$ is any non--constant holomorphic map $f:M \to N$, where $N$ a complex manifold, such that the fibers of $f$ are $\mathcal{F}$--invariant.
\end{remark}

\begin{theorem}\label{principal}
Let $\mathcal{F}$ be a $1$--dimensional foliation on the complex manifold $M$. If $V$ is a finite dimensional linear system, such that $\EC(\mathcal{F},V)$ vanishes identically then there exists an dense open set $U$ where $\mathcal{F}_{\mid U}$ admits a first integral. Moreover, if $M$ is a projective variety then $\mathcal{F}$ admits a meromorphic first integral.
\end{theorem}
\noindent{\it Proof}: Suppose that $\dim_{\C}C = k$, set $s_1,\ldots,s_k$ a basis of $V$ and consider $N \subset M \times \CPP$ defined as follows
\[
N = \left\lbrace (p ; a_1:\ldots:a_k) \mid \nu\left(\sum_{i=1}^{k}a_i \cdot s_i,\mathcal{F},p\right) = \infty \right\rbrace \, .
\]  

To see that $N$ is a closed set consider, for each positive integer $j$, the sets $N_j \subset M \times \CPP$ defined by
\[
N_j = \left\lbrace (p ; a_1:\ldots:a_k) \mid \nu\left(\sum_{i=1}^{k}a_i \cdot s_i,\mathcal{F},p\right) \ge j \right\rbrace \, .
\]
In an open set $U$ that trivializes $\mathcal{F}$ and the line bundle that supports $V$ we can write $N_j$ as
\[
N_j = \left\lbrace (p ; a_1:\ldots:a_k)  \mid X^l\left(\sum_{i=1}^{k}a_i \cdot s_i\right)(p)=0, l=0,\ldots,j \right\rbrace \, .
\]
Hence $N_j$ is a closed set. Since 
\[
  N = \bigcap_{j=1}^{\infty} N_j \, ,
\]
it follows that $N$ is also a closed set.

Observe that every fiber of the natural projection $\pi : N \to M$ is a projective space linearly embedded on $\CPP$. It is well known that there exists a dense open set $U$ where the dimension of the fibers are the same, say $l$. Hence, we can define a natural holomorphic function from $U \subset M$ to the Grassmanian of $(l+1)$--planes on $\C^k$. This function, defined in a dense open subset $U$ of $M$ is a first integral for $\mathcal{F}_{\mid U}$.

When $M$ is a projective variety we can extend this function to all $M$, obtaining in such way a meromorphic first integral to $\mathcal{F}$.
\cqd

%-------------------------------------------------------------------------------------------
% Generalizando Darboux
%----------------------------------------------------{--1}----------------------------------
\section{A generalization of Darboux's Criteria}
\rm
In \cite{L}, Jean-Marie Lion proposed the following generalization to Darboux's criteria for the existence of meromorphic first integral.

\vskip 0.2cm
\noindent {\bf Theorem[Lion]} { \it Let $\mathcal{F}$ be a holomorphic foliation given by an  integrable $1$-form, $\omega$, 
defined in a nighbourhood of the closure of $U$, where $U$ is a bounded open set of $\mathbb{C}^n$. Suppose that exists infinitely many leaves $V_i$ of  $\mathcal{F}$ contained in different algebraic hypersurfaces ${Q_i=0}$ of 
the same degree $d$. Then every leaf of $\mathcal{F}$ is contained in an irreducible algebraic 
hypersurface of degree at most $d$, and $\omega$ admits a meromorphic first integral. 
If $0 \in (\cap_i V_i)$, then exists a meromorphic first integral such that the graph is 
a $n$-dimensional algebraic subset of $\mathbb{C}^n \times \mathbb{P}^1_{\mathbb{C}}$.}
\vskip 0.2cm

Here, using the ideas developed earlier in this paper, we obtain a further generalization of the first part of Lion's result. The main improvement is that we do not need to restrict ourselves to algebraic leaves and open sets of $\C^n$.  We work with sections of any linear system defined on a  complex manifold which, in principle, do not have to be compact.

\begin{theorem}\label{genera}
Let $\mathcal{F}$ be a holomorphic foliation, of arbitrary codimension,on a complex manifold $M$, and $V$ 
finite dimensional linear system on $M$. Suppose that there exists a collection $\lbrace L_i\rbrace_{i \in \Lambda}$, such that each leaf $L_i$ has flat contact with some element of $V$. If the analytic closure of 
\[
  \bigcup_{i \in \Lambda} L_i \, 
\]
is equal to $M$ then $\mathcal{F}$ admits a first integral and every leaf of $\mathcal{F}$ has flat contact with some element of $V$.  
\end{theorem}
\noindent{\it Proof}:Let $U$ be a Stein open set of $M$. If we consider the $\mathcal{O}_U$--module $\mathcal{X}(\mathcal{F})$ formed by all vector fields in $U$ tangent to $\mathcal{F}_{\mid U}$, then for every $Y \in \mathcal{X}(\mathcal{F})$ we have that $\mathcal{E}(Y,V)=0$. Hence, as in the proof Theorem \ref{principal}, we can consider for each $Y \in \mathcal{X}(\mathcal{F})$ a set $N_{U,Y} \subset U \times \CPP$ describing the elements of $V$ with flat contact with $Y$. 

Defining $N \subset M \times \CPP$ locally as 
\[
  N_{\mid U \times \CPP} = \bigcap_{Y \in \mathcal{X}(\mathcal{F})} N_{U,Y} \, ,
\]
we obtain a closed analytic subset of $ M \times \CPP$ whose fibers under the natural projection to $M$ are projective spaces linearly embedded in $\CPP$. And as in Theorem \ref{principal} the result follows.
\cqd

\begin{corollary}
Let $\mathcal{F}$ be a codimension one holomorphic foliation on a complex manifold $M$, and $V$ 
finite dimensional linear system on $M$. Suppose that there exists an infinite collection $\lbrace L_i\rbrace_{i \in \Lambda}$, such that each leaf $L_i$ has flat contact with some element of $V$. If there exists an relatively compact open set $U \subset M$  such that $L_i \cap U \neq \emptyset$ for every $i \in \Lambda$ then $\mathcal{F}$ admits a first integral.
\end{corollary}

\section{Foliations with all leaves algebraic}

Gomez--Mont, in \cite{GM}, proved that a singular holomorphic foliation of codimension $q$ with all leaves algebraic admits a first integral, whose generic fiber has  codimension $q$. His proves uses Grothendieck's deep Theorem asserting the existence of the Hilbert scheme. Here, we propose a different proof using the ideas developed earlier in this paper.

\vskip 0.2cm
\noindent {\bf Theorem[Gomez--Mont]}
{\it Let $\mathcal{F}$ be a holomorphic foliation (with singularities) of codimension $q$ in the projective integral variety $M$, and assume that every leaf $L$ of $\mathcal{F}$ is a quasiprojective subvariety of $M$; then there is a projective integral variety $V$ of dimension $q$ and a rational map $f : M \to V$ such that the closure of a general $f$--fibre is the closure of a leaf of $\mathcal{F}$. }
\vskip 0.2cm
\noindent{\it Proof}: Suppose, without loss of generality, that $M$ is a subvariety of $\CPPP$ and take $\mathcal{L}$ to be the restriction of the hyperplane bundle $\mathcal{O}_{\mathbb{P}^n}(1)$ to $M$. 

Denote by $L_p$ the leaf through $p$ and $\rho(L_p)$  the least positive integer $l$ such that $L_p$ is an open subset of variety defined through sections of $\mathcal{L}^{\otimes l}$. Now, define
\[
  B_d = \lbrace p \in M \setminus {\rm Sing}(\mathcal{F}) \mid \rho(L_p) \le d  \rbrace \, .
\]

Since every leaf is algebraic, we have that 
\[
  M \setminus {\rm Sing}(\mathcal{F}) = \bigcup_{d\in \mathbb{N}} B_d \, ,
\]
and as consequence there exists $d_0$ such that $B_{d_0}$ has positive Lebesgue measure and consequently its analytic closure is equal to $M$.

Taking $V = H^0(M,\mathcal{L}^{\otimes d_0})$ we have from Theorem \ref{genera} that there exists a projective variety $N$ and a morphism $g:M \to N$ such that the fiber over any $p \in N$ is invariant by the foliation. 

Taking a closer look at the proof of Theorem \ref{genera} one can see that the fiber over a generic $p \in B_{d_0}$ will have codimension $q$. Using Stein's factorization Theorem, see \cite{Fi},  one can assure the existence of the commutative diagram,
\begin{equation*}
\begin{CD}
M @>f>> \overline{N} \\
@V{\rm Id}VV @V \pi VV \\
M @>g>> N
\end{CD} \, 
\end{equation*} 
where $\overline{N}$ is a projective  variety, $\pi$ is a  morphism and $f$ is a morphism whose generic fiber is irreducible. Hence $f$ satisfies the assertions of the Theorem.
\cqd

%-------------------------------------------------------------------------------------------
% Outras Aplicacoes
%-------------------------------------------------------------------------------------------

\section{Final Remarks}
The applications of the extactic curves and extactic divisors, certainly are not exhausted in this work. For example, in \cite{P} the extactic curves are used to introduce the  notion of algebraic multiplicity of an algebraic curve invariant by a polynomial vector field on $\C^2$. There, the relation of this algebraic multiplicity with the existence of exponential cofactors are explored to enrich Darboux's theory of integrability.

Another possible application is the study of automorphism groups of holomorphic foliations. The author plan to pursue this topic in a future paper.

%-------------------------------------------------------------------------------------------
%-------------------------------------------------------------------------------------------
%-------------------------------------------------------------------------------------------

%-------------------------------------------------------------------------------------------
% Referências
%-------------------------------------------------------------------------------------------

\end{document}